# ANALYTIC PROPERTIES OF MARKOV SEMIGROUP GENERATED BY STOCHASTIC DIFFERENTIAL EQUATIONS DRIVEN BY LÉVY PROCESSES

PANI W. FERNANDO, ERIKA HAUSENBLAS, AND PAUL RAZAFIMANDIMBY

ABSTRACT. We consider the stochastic differential equation (SDE) of the form
$$\begin{cases} dX^x(t) &= \sigma(X(t-))dL(t) \\ X^x(0) &= x, \quad x \in \mathbb{R}^d, \end{cases}$$
where $\sigma : \mathbb{R}^d \to \mathbb{R}^d$ is globally Lipschitz continuous and $L = \{L(t) : t \geq 0\}$ is a Lévy process. Under this condition on $\sigma$ it is well known that the above problem has a unique solution $X$. Let $(\mathcal{P}_t)_{t \geq 0}$ be the Markovian semigroup associated to $X$ defined by $(\mathcal{P}_t f)(x) := \mathbb{E}[f(X^x(t))]$, $t \geq 0$, $x \in \mathbb{R}^d$, $f \in \mathcal{B}_b(\mathbb{R}^d)$. Let $B$ be a pseudo–differential operator characterized by its symbol $q$. Fix $\rho \in \mathbb{R}$. In this article we investigate under which conditions on $\sigma$, $L$ and $q$ there exist two constants $\gamma > 0$ and $C > 0$ such that
$$|B\mathcal{P}_t u|_{H_2^\rho} \leq C\, t^{-\gamma}\, |u|_{H_2^\rho}, \quad \forall u \in H_2^\rho(\mathbb{R}^d), t > 0.$$

## 1. INTRODUCTION

The Blumenthal–Getoor index was first introduced in [6] in order to analyze the Hölder continuity of the sample paths, the $r$–variation, $r \in (0,2]$ and the Hausdorff-dimension of the paths of Lévy processes. Straightforward calculations give that the Blumenthal–Getoor index of an $r$–stable process is $r$. Lévy processes with Blumenthal–Getoor index less than 1 (resp. greater than 1) have paths of finite variation (resp. infinite variation). The Brownian motion has finite 2-variation. By using Hoh's symbol, Schilling introduced in [24], see also [4], a generalized Blumenthal–Getoor index which enabled him to characterize the Hölder continuity of the samples paths of a stochastic process. In [25], Schilling and Schnurr described the long term behavior in terms of this generalized index, for more details see also [4]. In [11] Glau gives a classification of Lévy processes via their symbols. To be more precise, Glau defines the Sobolev index of a Lévy process by a certain growth condition of its symbol.

In the present paper we investigate analytic properties of the Markovian semigroup generated by a SDE driven by a Lévy process. The main result in this article is Theorem 2.1 which is important, for instance, in nonlinear

Pani W. Fernando and Paul Razafimandimby were supported by the Austrian Science Foundations, Project number P 23591.

filtering with Lévy noise where one has to analyze the Zakai equation with jumps (see [10]). The leading operator of the Zakai equation is a pseudo–differential operator defined by the Hoh symbol of the driving noise in the state process. Thus, the uniqueness of the mild solution of the Zakai equation and its regularity depends very much on the estimate we obtain in Theorem 2.1.

Let us start our analysis by recalling the following definition (see [23]).

**Definition 1.1.** *A stochastic process $L = \{L(t) : t \geq 0\}$ with values in $\mathbb{R}^d$ is a called a d-dimensional Lévy process if the following conditions are satisfied.*

(1) *For any choice of $n \geq 1$ and $0 \leq t_0 < t_1 < \ldots < t_n$, the random variables $L(t_0), L(t_1) - L(t_0), L(t_2) - L(t_1), \ldots, L(t_n) - L(t_{n-1})$ are mutually independent (independent increment property).*
(2) $L(0) = 0$ *a.s.*
(3) *The distribution of $L(t+s) - L(s)$ does not depend on $s$ (stationary increment property).*
(4) $L$ *is stochastically continuous.*
(5) *There is $\Omega_0 \in \mathcal{F}$ with $\mathrm{P}[\Omega_0] = 1$ such that, for every $\omega \in \Omega_0$, the trajectory $[0, \infty) \ni t \mapsto L(t, \omega)$ is right-continuous and has left limits.*

Let $L = \{L(t) : t \geq 0\}$ be a $d$–dimensional Lévy process. We consider the stochastic differential equations of the form

$$(1.1) \quad \begin{cases} dX^x(t) &= b(X^x(t-))\,dt + \sigma(X^x(t-))\,dL(t) \\ X^x(0) &= x, \quad x \in \mathbb{R}^d, \end{cases}$$

where $b : \mathbb{R}^d \to \mathbb{R}^d$ and $\sigma : \mathbb{R}^d \to \mathbb{R}^d \times \mathbb{R}^d$ are functions satisfying the following conditions.

**Hypothesis 1.** *Let $k \geq \frac{d}{2}$. We assume that $b \in C^k(\mathbb{R}^d; \mathbb{R}^d)$ and $\sigma \in C^k(\mathbb{R}^d; \mathbb{R}^d \times \mathbb{R}^d)$ and they with their derivatives are bounded. In particular, we suppose that $\sigma$ is bounded from below and above.*

Under Hypothesis 1, the existence and uniqueness of a solution to equation (1.1) is well established, see for e.g. [2, p. 367, Theorem 6.2.3]. In addition, $X^x$ has $\mathbb{P}$–a.s. càdlàg trajectories.

Let $(\mathcal{P}_t)_{t \geq 0}$ be the Markovian semigroup associated to $X$ defined by

$$(\mathcal{P}_t f)(x) := \mathbb{E}\left[f(X^x(t))\right], \ \ t \geq 0, \ \ x \in \mathbb{R}^d, \ \ f \in \mathcal{B}_b(\mathbb{R}^d).$$

We have the following results.

**Proposition 1.1.** *If $\sigma$ and $b$ satisfies Hypothesis 1, then*

(1) *for any $t, s \geq 0$ we have $\mathcal{P}_t \circ \mathcal{P}_s = \mathcal{P}_{t+s}$.*
(2) *the semigroup $(\mathcal{P}_t)_{t \geq 0}$ is Feller on $C_b(\mathbb{R}^d)$ and $C_0(\mathbb{R}^d)$.*

*Proof.* Since, by Hypothesis 1, $\sigma$ and $b$ are globally Lipschitz, the proof of item (1) is quite standard and can be found, for instance, in [2, Section 6.4.2]. Owing to Hypothesis 1 again, the $C_0$-Feller of $(\mathcal{P}_t)_{t \geq 0}$ follows from

[2, Theorem 6.7.2]. For the $C_b$-Feller property we refer, for instance, to [2, Note 3]. □

The infinitesimal generator of $(\mathcal{P}_t)_{t\geq 0}$ is given by

$$Au(x) = -\int_{\mathbb{R}^d} e^{ix^T\xi}\psi(x,\xi)\hat{u}(\xi)\,d\xi \quad u \in C_b^\infty(\mathbb{R}^d),$$

where the symbol $\psi$ is defined by

$$\psi(x,\xi) := -\lim_{t\downarrow 0}\frac{1}{t}\mathbb{E}\left[e^{i(X^x(t)-x)^T\xi} - 1\right], \quad x \in \mathbb{R}^d.$$

In case $L$ is a $d$–dimensional Brownian motion and $\sigma \in C_b^\infty(\mathbb{R}^d;\mathbb{R}^d\times\mathbb{R}^d)$ is bounded from below and above, $A$ is a second order partial differential operator on $L^2(\mathbb{R}^d)$ with domain $D(A) = H_2^2(\mathbb{R}^d)$. Moreover, $(\mathcal{P}_t)_{t\geq 0}$ is an analytic semigroup on $L^2(\mathbb{R}^d)$ and the following inequality holds for $B = \nabla$

$$(1.2) \qquad |B\mathcal{P}_t x|_{L^2} \leq \frac{1}{\sqrt{t}}|x|_{L^2}, \quad x \in L^2(\mathbb{R}^d), \quad t > 0.$$

Let $L$ be a pure jump Lévy process and $B$ a pseudo–differential operator induced by a symbol. The purpose of this article is to investigate under which additional conditions the estimate (1.2) holds.

**Notation 1.1.** *For any nonnegative integers $\alpha_1, \alpha_2, \ldots, \alpha_d$ we set $\alpha = (\alpha_1, \ldots, \alpha_d)$ and $|\alpha| = \sum_{j=1}^d \alpha_j$. Moreover for a function $f : \mathbb{R}^d \to \mathbb{C}$ we write $\partial_x^\alpha f(x)$ for*

$$\frac{\partial^\alpha}{\partial x_1 \partial x_2 \cdots \partial x_d}f(x).$$

*For any $\rho \geq 0$ we define the function $\langle \cdot \rangle : \mathbb{R} \ni \xi \mapsto \langle \xi \rangle^\rho := (1+|\xi|^2)^{\frac{\rho}{2}} \in \mathbb{R}$. The following inequality, also called the Peetre inequality is used in several places: for any $s \in \mathbb{R}$ there exists a constant $c_s > 0$ such that*

$$\langle x + y \rangle^s \leq c_s \langle x \rangle^s \langle y \rangle^{|s|}, \quad x, y \in \mathbb{R}^d.$$

Let $U \subset \mathbb{R}^d$ be a non–empty set and $f, g : U \to [0, \infty)$. We set $f(x) \lesssim g(x)$, $x \in U$, iff there exists a constant $C > 0$ such that $f(x) \leq Cg(x)$ for all $x \in U$. Moreover, if $f$ and $g$ depend on a further variable $z \in Z$, the statement for all $z \in Z$, $f(x,z) \lesssim g(x,z)$, $x \in U$ means that for every $z \in Z$ there exists a real number $C_z > 0$ such that $f(x,z) \leq C_z g(x,z)$ for every $x \in U$. Also we set $f(x) \asymp g(x)$, $x \in U$, iff $f(x) \lesssim g(x)$ and $g(x) \lesssim f(x)$ for all $x \in U$. Finally, we say $f(x) \gtrsim g(x)$, $x \in U$, iff $g(x) \lesssim f(x)$, $x \in U$. Similarly as above, one may handle the case if the functions depend on a further variable.

Let $\mathcal{S}(\mathbb{R}^d)$ be the Schwartz space of functions $C^\infty(\mathbb{R}^d)$ where all derivatives decreases faster than any power of $|x|$ as $|x|$ approaches infinity. Let $\mathcal{S}'(\mathbb{R}^d)$ be the dual of $\mathcal{S}(\mathbb{R}^d)$. For any pair of functions $f, g \in \mathcal{S}(\mathbb{R}^d)$, we define $(f,g)$ by

$$(f,g) := \int_{\mathbb{R}^d} f(x)\bar{g}(x)\,dx.$$

Throughout this paper, we define the Fourier transform $\mathcal{F}$ as follows:

$$\mathcal{F}f(\xi) = \hat{f}(\xi) = (2\pi)^{-d}\int_{\mathbb{R}^d} e^{i\xi^T x}f(x)\,dx, \quad f \in L^2(\mathbb{R}^d).$$

Its inverse is $\mathcal{F}^{-1}$ and is defined by

$$\mathcal{F}^{-1}f(x) = \check{f}(x) = \int_{\mathbb{R}^d} e^{i\xi^T x}f(\xi)\,d\xi, \quad f \in L^2(\mathbb{R}^d).$$

Let $1 \leq p < \infty$ and $s \in \mathbb{R}$, then $H_p^s(\mathbb{R}^d)$ denotes the Bessel Potential spaces or Sobolev spaces of fractional order defined by (see [29, Chapter 3.2])

$$H_p^s(\mathbb{R}^d) := \left\{ f \in L^p(\mathbb{R}^d) : |f|_{H_m^p} = |\mathcal{F}^{-1}\langle\xi\rangle^s(\mathcal{F}f)|_{L^p} < \infty \right\}.$$

Let $(U, d)$ be a metric space. By $C(U)$ (resp. $\mathcal{B}_b(U)$ ) we denote the set of all complex valued continuous (resp. measurable and bounded) functions $f : U \to \mathbb{C}$. Furthermore, for $m \in \mathbb{N}$ we define

$$C^{(m)}(U) = \{f \in C(U) : \partial_x^\alpha f \in C(U) \text{ for all } |\alpha| \leq m\}.$$

In addition, let $C_b(\mathbb{R}^d) = \{f \in C(\mathbb{R}^d) : \sup_{x \in \mathbb{R}^d}|f(x)| < \infty\}$ and let $C_b^k(\mathbb{R}^d) = \{f \in C(U) : \partial_x^\alpha f \in C_b(\mathbb{R}^d) \text{ for all } |\alpha| \leq m\}$. The symbol $C_0(U)$ denotes the set of continuous functions $f : U \to \mathbb{C}$ such that

$$\lim_{|y|\to\infty} f(y) = 0.$$

Finally, we denote by $\mathbb{N}$ the set of positive integers and $\mathbb{N}_0 = \mathbb{N} \cup \{0\}$.

## 2. Hoh's Symbols associated to Lévy processes

Throughout the remaining article, let $L = \{L(t) : t \geq 0\}$ be a family of $d$–dimensional Lévy processes and let us denote by $L^x = L + x$, $x \in \mathbb{R}^d$. Then $L^x$ generates a Markovian semigroup $(\mathcal{T}_t)_{t\geq 0}$ on $C_b(\mathbb{R}^d)$ by

$$\mathcal{T}_t f(x) := \mathbb{E}f(L^x(t)), \quad f \in C_b(\mathbb{R}^d).$$

Let $A$ be the infinitesimal generator of $(\mathcal{T}_t)_{t\geq 0}$ defined by

$$(2.1) \qquad Af := \lim_{h\to 0}\frac{1}{h}\left(\mathcal{T}_h - \mathcal{T}_0\right)f, \quad f \in C_b^{(2)}(\mathbb{R}^d).$$

An alternative way of defining $A$ makes use of the Lévy symbols. In particular, let $\psi : \mathbb{R}^d \to \mathbb{C}$ be defined by

$$\psi(\xi) = -\lim_{t\downarrow 0}\frac{1}{t}\ln(\mathbb{E}e^{i\xi^T(L^x(t)-x)}), \quad \xi \in \mathbb{R}^d.$$

Then, there exist a vector $l \in \mathbb{R}^d$, a positive definite matrix $Q$ and a Lévy measure $\nu$ such that $\psi$ can be written in the following form

$$\psi(\xi) = il^T\xi - \frac{1}{2}\xi^T Q\xi + \int_{\mathbb{R}^d\setminus\{0\}} \left(e^{i\xi^T z} - 1 - i\xi^T z 1_{|z|\leq 1}\right)\nu(dz), \quad \xi \in \mathbb{R}^d.$$

Here, $\psi$ is called the Lévy symbol of the Lévy process $L$ (see [23]). It can be shown that if $L$ is a Lévy process with symbol $\psi$, then the infinitesimal generator defined by (2.1) can also be written as (see e.g. [2, 17, 4])

(2.2)
$$A f(x) = -\int_{\mathbb{R}^d} e^{i\xi^T x}\psi(\xi)\hat{f}(\xi)\,d\xi, \quad x \in \mathbb{R}^d,\, f \in \mathcal{S}(\mathbb{R}^d).$$

If $\nu$ has bounded second moments, the operator $A$ is well defined on $C_b^2(\mathbb{R}^d)$, has values in $\mathcal{B}_b(\mathbb{R}^d)$ and satisfies the positive maximum principle (see e.g. [17, Theorem 4.5.13 ]). Therefore, $A$ generates a Feller semigroup on $C_b^\infty(\mathbb{R}^d)$ and a sub–Markovian semigroup on $L^2(\mathbb{R}^d)$ (see e.g. [18, Theorem 2.6.9 and Theorem 2.6.10]).

A further important property of a Lévy symbol is, that it is negative definite.

**Definition 2.1.** *A function $\psi : \mathbb{R}^d \to \mathbb{C}$ is called negative definite function iff for all $m \in \mathbb{N}$ and all $m$–tuples $(\xi_1, \ldots, \xi_m)$, $\xi^j \in \mathbb{R}^d$, $1 \leq j \leq m$, the matrix*

$$(\psi(\xi_j) + \overline{\psi(\xi_j)} - \psi(\xi^k - \xi^j))_{k,j=1,\cdots m}$$

*is positive Hermitian, i.e. for all $c_1, \ldots, c_m \in \mathbb{C}$*

$$\sum_{j,k=1}^m (\psi(\xi_j) + \overline{\psi(\xi_j)} - \psi(\xi^k - \xi^j))c_k\bar{c}_j \geq 0.$$

If a negative definite function $\psi$ is continuous, then

$$|\psi(\xi)| \lesssim \langle|\xi|^2\rangle, \quad \xi \in \mathbb{R}^d.$$

Moreover, every real-valued continuous negative definite function $\psi$ has a representation

$$\psi(\xi) = c + q(\xi) + \int_{\mathbb{R}^d\setminus\{0\}} (1 - \cos(\xi^T y))\nu(dy),$$

where $c \geq 0$ is a constant, $q \geq 0$ is a quadratic form and $\nu$ is a symmetric Borel measure on $\mathbb{R}^d \setminus \{0\}$ called Lévy measure having the property that

$$\int_{\mathbb{R}^d\setminus\{0\}} \frac{|y|^2}{1+|y|^2}\nu(dy) < \infty.$$

The Blumenthal–Getoor index of a Lévy process with Lévy measure $\nu$ is defined by the

$$\inf\{p : \int_{|x|\leq 1} |x|^p\nu(dx) < \infty\}.$$

The Blumenthal–Getoor index has been generalized in several direction, see e.g. [4, p. 124]. Here, in this article we modify also the Blumenthal–Getoor index, but similarly to the characterization of pseudo–differential operators.

**Definition 2.2.** *Let $L$ be a normalized Lévy process with symbol $\psi$ and $\psi \in C^k(\mathbb{R}^d \setminus \{0\})$ for some $k \in \mathbb{N}_0$. Then the Blumenthal–Getoor index of order $k$ is defined by*

$$s := \inf_{\substack{\lambda > 0 \\ |\alpha| \le k}} \left\{ \lambda : \lim_{|\xi| \to \infty} \frac{|\partial_\xi^\alpha \psi(\xi)|}{|\xi|^{\lambda - |\alpha|}} = 0 \right\}.$$

*Let*

$$s^+ := \inf_{\substack{\lambda > 0 \\ |\alpha| \le k}} \left\{ \lambda : \limsup_{|\xi| \to \infty} \frac{|\partial_\xi^\alpha \psi(\xi)|}{|\xi|^{\lambda - |\alpha|}} = 0 \right\}$$

*be the upper and*

$$s^- := \inf_{\substack{\lambda > 0 \\ |\alpha| \le k}} \left\{ \lambda : \liminf_{|\xi| \to \infty} \frac{|\partial_\xi^\alpha \psi(\xi)|}{|\xi|^{\lambda - |\alpha|}} = 0 \right\}$$

*the lower Blumenthal–Getoor index $s^+$ of order $k$. Here $\alpha$ denotes a multi-index. If $k = \infty$ then Blumenthal–Getoor index of infinity order is defined by*

$$s := \inf_{\substack{\lambda > 0 \\ \alpha \text{ is a multi-index}}} \left\{ \lambda : \lim_{|\xi| \to \infty} \frac{|\partial_\xi^\alpha \psi(\xi)|}{|\xi|^{\lambda - |\alpha|}} = 0 \right\}.$$

**Remark 2.1.** *The Blumenthal–Getoor index of order infinity is defined for the sake of completeness. In this paper, we are interested in Lévy processes with Blumenthal–Getoor index of finite order $k$.*

In order to define the resolvent of an operator $A$ associated to a symbol $\psi$, we need to characterize the range of the symbol of $\psi$. Here, one has different possibilities at one's disposal, depending on whether one comes from analysis or probability theory. Thus, one can introduce the sector condition or the type of a symbol. Both definitions describe the range of a symbol. The sector condition reads as follows.

**Definition 2.3.** *(compare [3, p. 116]) Let $L$ be a Lévy process with symbol $\psi$. We say the symbol $\psi$ satisfies the sector condition, if there exists a $\kappa > 0$ such that*

$$|\Im(\psi(\xi))| \le \kappa \, \Re(\psi(\xi)), \quad \forall x, \xi \in \mathbb{R}^d.$$

For $\delta \in [0, \pi]$ we define $\Sigma_\delta := \{ z \in \mathbb{C} \setminus \{0\} : |\arg(z)| < \delta \}$. Now, the type of a symbol is given as follows.

**Definition 2.4.** *Let $L$ be a Lévy process with symbol $\psi$. We say the symbol $\psi$ is of type $(\omega, \theta)$, $\omega \in \mathbb{R}$, $\theta \in (0, \frac{\pi}{2})$, iff*

$$-\overline{\psi(\mathbb{R}^d)} \subset \mathbb{C} \setminus \omega + \Sigma_{\theta + \frac{\pi}{2}}.$$

**Remark 2.2.** *In general, one uses the essential range of a symbol $\psi$ to characterize its type. But in our case $\psi$ is continuous so we work with the closure of the range of $\psi$.*

**Remark 2.3.** *If a symbol $\psi$ is of type $(0,\theta)$, then it satisfies the sector condition with $\kappa = \tan(\theta)$ and vice versa.*

**Remark 2.4.** *Let $R(\lambda : A) = [\lambda + A]^{-1}$ be the resolvent of the operator $A$ and let $M_\psi$ be the multiplication operator induced by the function $\psi$ (see e.g. [7, Definition 4.1, p. 25]). First, note that by the spectral Theorem (see e.g. [7, Theorem 4.9, p. 30]), $A$ is unitary equivalent on $L^2(\mathbb{R}^d)$ to the multiplication operator $M_\psi$. To be more precise, there exists a unitary operator, in our case the Fourier transform $\mathcal{F} : L^2(\mathbb{R}^d) \to L^2(\mathbb{R}^d)$, such that $A = \mathcal{F}^{-1}M_\psi(\mathcal{F}f)$. Similarly, if $\lambda \in \mathbb{C}$ belongs to the spectrum of $M_\psi$, then $(\lambda I - M_\psi)$ is not invertible. It follows that the operator $\lambda I - A = \mathcal{F}^{-1}(\lambda - M_\psi)\mathcal{F}$ is not invertible, and therefore $\lambda$ belongs also to the spectrum of $A$. Note, by Proposition 4.2 [7, p. 25] it follows that we have for the spectrum $\sigma(M_\psi)$ of the operator $M_\psi$*

$$\sigma(M_\psi) = \overline{\psi(\mathbb{R}^d)}.$$

*In addition, by Theorem 1.4.2 of [13] it follows for $\lambda \in \mathbb{C} \setminus \overline{\psi(\mathbb{R}^d)}$*

$$\|R(\lambda, M_\psi)\|_{L(L^2)} \leq \frac{1}{dist(\overline{\psi(\mathbb{R}^d)}, \lambda)},$$

*which implies*

$$\|R(\lambda, A)\|_{L(L^2)} \leq \frac{1}{dist(\overline{\psi(\mathbb{R}^d)}, \lambda)}$$

*In addition, taking into account that the multiplication operator $m(\xi) = \langle\xi\rangle^s$ is an isomorphism from $L^2(\mathbb{R}^d)$ to $H^s(\mathbb{R}^d)$ one has for any $s \in \mathbb{R}$*

$$\|R(\lambda, A)\|_{L(H^s)} \leq \frac{1}{dist(\overline{\psi(\mathbb{R}^d)}, \lambda)}.$$

The generalized Blumenthal–Getoor index of order 0 and the type of a symbol can be calculated in many cases. Here, we give some examples.

**Example 2.1.** *Let $\alpha \in (0,2)$ and $L$ be a symmetric $\alpha$–stable process without drift. The symbol $\psi$ of $L$ is given by*

$$\psi(\xi) = |\xi|^\alpha,$$

*the upper and lower index is $\alpha$, and $\psi$ is of type $(0,\delta)$ for any $\delta > 0$.*

**Example 2.2.** *Let $L$ be the Meixner process as described in [26] (see also [19, p. 136] or [3]). In particular, let $L$ be a real–valued Lévy process with symbol*

$$\psi_{m,\delta,a,b}(\xi) = -im\xi + 2\delta\left(\log\cosh\left(\frac{a\xi - ib}{2}\right) - \log\cos\left(\frac{b}{2}\right)\right), \quad \xi \in \mathbb{R},$$

where $m \in \mathbb{R}$, $\delta, a > 0$, $b \in (-\pi, \pi)$. Then the upper and lower index is 1 (see [19, p. 137-(3.226)]). Moreover, the symbol $\psi$ is of type $(\omega, \theta)$ with $\omega = 0$ and $\theta = \arctan(m/\delta a)$.

**Example 2.3.** *Let $L$ be the normal inverse Gaussian process as described in [5] (see also [19, p. 138]). In particular, let $L$ be a real–valued Lévy process with symbol*

$$\psi_{NIG}(\xi) = -im\xi + \delta\left(\sqrt{a^2 - (b + i\xi)} - \sqrt{a^2 - b^2}\right), \quad \xi \in \mathbb{R},$$

*where $m \in \mathbb{R}$, $\delta > 0$, $0 < |b| < a$. This process is comparable with the Cauchy process, but has finite expectation. Next, the upper and lower index is 1 (see [19, p. 137-(3.228)]). Moreover, for $m = 0$ the symbol $\psi$ is of type $(\omega, \theta)$ with $\omega = 0$ and $\theta < \pi$.*

**Example 2.4.** *It is interesting to note that, by subordination, one can produce sectorial symbols given a non-sectorial one. The typical example $\psi(\xi) = i\xi$ (non-sectorial) becomes sectorial if we use $(i\xi)^\alpha$, i.e. the Bernstein function $f(s) = s^\alpha$ with $\alpha \in (0, 1)$. There are many of Bernstein functions having this property, for references see e.g. [12, 16].*

Let $L = \{L(t) : t \geq 0\}$ be a $d$–dimensional Lévy process without any Gaussian component. We consider the stochastic differential equations of the form

$$(2.3) \quad \begin{cases} dX^x(t) &= \sigma(X^x(t-))dL(t) \\ X^x(0) &= x, \quad x \in \mathbb{R}^d, \end{cases}$$

where $\sigma : \mathbb{R}^d \to \mathbb{R}^d \times \mathbb{R}^d$ is map satisfying Hypothesis 1. Let $(\mathcal{P}_t)_{t \geq 0}$ be the associated Markovian semigroup of $X$ defined by

$$(2.4) \quad (\mathcal{P}_t f)(x) := \mathbb{E}\left[f(X^x(t))\right], \quad t \geq 0, \quad f \in C_b(\mathbb{R}^d).$$

Then, $(\mathcal{P}_t)_{t \geq 0}$ is a Feller semigroup and one can compute its infinitesimal generator. Again, one way of computing $A$ is done by Hoh's symbols (see [15]). In particular, one has

$$Au(x) = \int_{\mathbb{R}^d} e^{ix^T\xi} a(x, \xi) \hat{u}(\xi)\, d\xi \quad u \in C_b^\infty(\mathbb{R}^d),$$

where the symbol $p$ is defined by

$$a(x, \xi) := -\lim_{t \downarrow 0} \frac{1}{t} \mathbb{E}\left[e^{i(X^x(t)-x)^T\xi} - 1\right], \quad x \in \mathbb{R}^d.$$

Alternatively, we can give an explicit form of $a(\cdot, \cdot)$ in term of the Lévy symbol of the driving noise $L$. In fact, if $\psi$ is a Lévy symbol of the Lévy process $L = \{L(t) : t \geq 0\}$, then it is shown in [25, Theorem 3.1], that the symbol $a : \mathbb{R}^d \times \mathbb{R}^d \to \mathbb{C}$ is given by

$$(2.5) \quad a(x, \xi) = \psi(\sigma^T(x)\xi), \quad (x, \xi) \in \mathbb{R}^d \times \mathbb{R}^d.$$

Symbols also arise in the context of pseudo–differential operators, whereas the term symbol is defined in the following way (in Appendix B we give a short summary of some definition and theorems that we need for our proof).

**Definition 2.5.** *(compare* [30, p.28, Def. 4.1]*,* [21, Def. 1.1.1, p. 19]*,* [27, p.3 Def.1.1]*) Let $X \subset \mathbb{R}^d$, $m \in \mathbb{R}$, and $\rho, \delta$ two real numbers such that $0 \le \rho \le 1$ and $0 \le \delta \le 1$. Let $S^m_{\rho,\delta}(X, \mathbb{R}^d)$ be the set of all functions $a : X \times \mathbb{R}^d \to \mathbb{C}$, where*

- *$a$ is infinitely often differentiable;*
- *for any two multi-indices $\alpha$ and $\beta$ there exists a positive constant $C_{\alpha,\beta} > 0$ depending only on $\alpha$ and $\beta$ such that*

$$\left|\partial_x^\alpha \partial_\xi^\beta a(x,\xi)\right| \le C_{\alpha,\beta} \langle |\xi| + |x| \rangle^{m-\rho|\beta|+\delta|\alpha|}, \quad x \in X, \xi \in \mathbb{R}^d.$$

The symbolic calculus for pseudo–differential operators is well established, see [21, 27, 30]. However, in case one considers symbols associated to the solution to stochastic differential equations driven by Lévy processes, the derivatives of the symbol will not be necessarily continuous at the origin, i.e. at $\{0\}$. The behavior of $\xi$ at the origin corresponds to the perturbation of the solution $X$ of equation (2.3) by the large jumps of the Lévy process $L$. To illustrate this fact, let us assume that the symbol $a(x,\xi) = a(\xi)$ is independent from $x$ and positive definite. Then the symbol corresponds to a Lévy process. Now, let us assume that the Lévy process has a symmetric Lévy measure $\nu$ such that for all $\ell \ge 2$ the moments $\int_{\mathbb{R}^d \setminus \{0\}} |y|^\ell \nu(dy)$ are bounded. Then by [15, Proposition 2; p.793] the Lévy symbol $a$ is infinitely often differentiable and one has

$$\left|\partial_\xi^\alpha a(\xi)\right| \lesssim \begin{cases} a(\xi) & \text{if } \alpha = 0, \\ |a(\xi)|^{\frac{1}{2}} & \text{if } |\alpha| = 1, \alpha \in \{1,\ldots,d\} \\ 1 & \text{if } |\alpha| \ge 2, \end{cases}$$

which means that if all moments of the Lévy measure are bounded, i.e. the moments of the large jumps are bounded, then the Lévy symbol will be infinitely often differentiable at the origin. Now, let us assume that the Lévy process is symmetric and $r$–stable with $r < 2$. It is well known that the Lévy symbol is $|\xi|^r$ and the large jumps have only bounded moments up to order $\ell$ with $\ell < r$. In case $r < 1$ the Lévy symbol is only once continuously differentiable at the origin and in case $r > 1$, twice continuous differentiable at the origin.

Now, one may ask the question: does the non differentiability at the origin have any effect on the smoothing property of the corresponding Markovian semigroup $(\mathcal{P}_t)_{t \ge 0}$. Again, let us assume that the Lévy process is symmetric and $r$–stable with $r < 2$. Then, the infinitesimal generator of the corresponding Markovian semigroup $(\mathcal{P}_t)_{t \ge 0}$ is $-(-\Delta)^{\frac{r}{2}}$. However, it is well known that

$$|\mathcal{P}_t u|_{H_2^r} \le \frac{c}{t} |u|_{L^2},$$

which implies that the discontinuity of the derivatives at the origin will not have an effect on the smoothing property of the corresponding semigroup $(\mathcal{P}_t)_{t\geq 0}$. However, the discontinuity at the origin has to be taken into account and we will relax the definition of symbols slightly and define a wider class of Hoh symbols.

**Definition 2.6.** Let $m \in \mathbb{R}$, and $\rho, \delta$ two real numbers such that $0 \leq \rho \leq 1$ and $0 \leq \delta \leq 1$ and $k \in \mathbb{N}_0$. Let $\mathcal{S}^m_{h k; \rho, \delta}(\mathbb{R}^d, \mathbb{R}^d)$ be the set of all functions $a : \mathbb{R}^d \times (\mathbb{R}^d \setminus \{0\}) \to \mathbb{C}$, where

- for all multi-indices $\alpha$ and $\beta$ with $|\alpha|, |\beta| \leq k$ one has $\partial_\xi^\alpha \partial_x^\beta a \in C(\mathbb{R}^d \times \mathbb{R}^d \setminus \{0\})$ and

$$\left| \partial_\xi^\alpha \partial_x^\beta a(x, \xi) \right| \lesssim |\xi|^{-|\alpha|}, \quad |\xi| \leq 1;$$

- for any two multi-indices $\alpha$ and $\beta$ with $|\alpha|, |\beta| \leq k$, there exists a positive constant $C_{\alpha, \beta} > 0$ depending only on $\alpha$ and $\beta$ such that

$$\left| \partial_x^\alpha \partial_\xi^\beta a(x, \xi) \right| \leq C_{\alpha, \beta} \langle |\xi| + |x| \rangle^{m - \rho|\beta| + \delta|\alpha|}, \quad x, \xi \in \mathbb{R}^d, |\xi| \geq 1.$$

**Remark 2.5.** Let us assume that a Lévy symbol $\psi$ has a generalized Blumenthal–Getoor index $s$ of order $k \geq 1$ and $\sigma \asymp 1$, $\sigma \in C_b^k(\mathbb{R}^d)$. Direct computation shows that the symbol $a(x, \xi) := \psi(\sigma^T(x)\xi)$ belongs to $\mathcal{S}^s_{h k; 1, 0}(\mathbb{R}^d, \mathbb{R}^d)$.

We also make the following observation.

**Remark 2.6.** Hypothesis 1 can be reformulated as follows. We assume that $b \in C_b^k(\mathbb{R}^d)$ and $\sigma \in C_b^k(\mathbb{R}^d; \mathbb{R}^d \times \mathbb{R}^d)$ such that for all multi-indices $|\alpha| \leq k$ one has $|\partial_x^\alpha \sigma(x)| \asymp 1$ and $|\partial_x^\alpha b(x)| \lesssim 1$, $x \in \mathbb{R}^d$.

Now, we can formulate the following Theorem.

**Theorem 2.1.** Let $k \geq \frac{d}{2}$ be a real number, and $b : \mathbb{R}^d \to \mathbb{R}^d$ and $\sigma : \mathbb{R}^d \to \mathbb{R}^d \times \mathbb{R}^d$ be two maps satisfying Hypothesis 1. Let $q$ and $\psi$ be two generalized Hoh symbols, $k$ times differentiable on $\mathbb{R}^d \setminus \{0\}$. Let $X$ be the unique solution to (1.1), where $L$ is a Lévy process with symbol $\psi$. Let $B : \mathcal{S}(\mathbb{R}^d) \to \mathcal{S}'(\mathbb{R}^d)$ be an operator defined by

$$(Bu)(x) := -\lim_{t \downarrow 0} \int_{\mathbb{R}^d} e^{ix^T \xi} q(b^T(x)\xi) \hat{u}(\xi) \, d\xi, \quad u \in \mathcal{S}(\mathbb{R}^d).$$

If $q$ and $\psi$ are $k$–times differentiable on $\mathbb{R}^d$, then we do not impose any further condition on $\psi$ and $q$. If $\psi$ and $q$ are only $k$–times differentiable on $\mathbb{R}^d \setminus \{0\}$, we assume that $\psi(0) = 0$, and $q(0) = 0$, and there exist two constants $\gamma_\psi > 0$ and $\gamma_q > 0$ for which one has

$$\sup_{|\xi| \leq 1, \xi \neq 0} \left| |\xi|^{|\alpha| - \gamma_\psi} \partial_\xi^\alpha \psi(\xi) \right| < \infty,$$

and

$$\sup_{|\xi| \leq 1, \xi \neq 0} \left| |\xi|^{|\alpha| - \gamma_q} \partial_\xi^\alpha q(\xi) \right| < \infty,$$

*for all multi-indices $\alpha$ with $|\alpha| \leq k$.*

If $\psi$ is of type $(\omega, \theta)$ and has generalized Blumenthal–Getoor index $s_1$ and $q$ has upper Blumenthal–Getoor index less or equal to $s_2$ with $s_2 < s_1$, then for any $\rho \in \mathbb{R}$ there exists a constant $C > 0$ such that

$$(2.6) \quad |B\,\mathcal{P}_t u|_{H_2^\rho} \leq \frac{C}{\sin\theta}\, t^{-\frac{s_2}{s_1}} |u|_{H_2^\rho}, \quad u \in H_2^\rho(\mathbb{R}^d),\, t > 0.$$

**Remark 2.7.** *In the following we will write $a(x, \xi)$ to denote the symbol and $a(x, D)$ to denote the corresponding pseudodifferential operator given by*

$$a(x, D)u(x) := \int_{\mathbb{R}^d} e^{ix^T \xi} a(x, \xi) \hat{u}(\xi)\, d\xi, \quad u \in \mathcal{S}(\mathbb{R}^d).$$

**Remark 2.8.** *The same holds if one of the operators is replaced by its adjoint. In fact, for a symbol $a \in S_{\delta,\rho}^m(\mathbb{R}^d, \mathbb{R}^d)$ we know by Theorem B.2, that there exists a symbol $a^*(x, \xi) \in S_{\delta,\rho}^m(\mathbb{R}^d, \mathbb{R}^d)$ such that the adjoint operator $a^*(x, D)$ of $a(x, D)$ is described by the symbol $a^*(x, \xi)$. In addition, by Remark B.2, we know that if $a \in \operatorname{Hyp}_{\delta,\rho}^{m,m_0}(\mathbb{R}^d, \mathbb{R}^d)^1$, then $a^*(x, \xi) \in \operatorname{Hyp}_{\delta,\rho}^{m,m_0}(\mathbb{R}^d, \mathbb{R}^d)$. Thus, if the symbol $\psi(\sigma^T(x)\xi)$ and $q(b^T(x)\xi)$ satisfy the assumption of Theorem 2.1, then $\psi^*(\sigma^T(x)\xi)$ and $q^*(b^T(x)\xi)$, respectively, satisfy also the assumption of Theorem 2.1.*

*Proof of Theorem 2.1.* First, note that we use within the proof the notation introduced in appendix B.

Let $a(x, \xi) := \psi(\sigma(x)^T \xi)$, $A = a(x, D)$, $\operatorname{Dom}(A) = \{f \in L^2(\mathbb{R}^d) : Af \in L^2(\mathbb{R}^d)\}$. Let us denote by $R(\lambda : A)$ the resolvent of the operator $A$. The proof of Theorem 2.1 relies on Proposition A.1 and the fact that there exists a constant $C > 0$ such that for $\varepsilon = s_2/s_1$

$$(2.7) \quad |B\,R(\lambda : A)u|_{H_2^\rho} \leq C\, |\lambda|^{\varepsilon - 1}\, |u|_{H_2^\rho}, \quad \forall u \in H_2^\rho(\mathbb{R}^d).$$

In particular, if (2.7) is true, then an application of Proposition A.1 gives the assertion. Hence, we are going to show that under the conditions of Theorem 2.1 estimate (2.7) holds.

In the first step, we assume that $\psi$, $q$, $\sigma$ and $b$ are Schwartz functions defined on $\mathbb{R}^d$, i.e. $\psi, q, \sigma, b \in \mathcal{S}(\mathbb{R}^d)$. In the second step we replace $\psi$ and $q$ by sequences $\{\psi_n : n \in \mathbb{N}\} \subset \mathcal{S}(\mathbb{R}^d)$ and $\{q_n : n \in \mathbb{N}\} \subset \mathcal{S}(\mathbb{R}^d)$ converging respectively to $\psi$ and $q$. We will also replace $\sigma$ and $b$ by sequences $\{\sigma_n : n \in \mathbb{N}\} \subset \mathcal{S}(\mathbb{R}^d)$ and $\{b_n : n \in \mathbb{N}\} \subset \mathcal{S}(\mathbb{R}^d)$ converging respectively to $\sigma$ and $b$, in an appropriate sense that will be made precise later.

**Step 1:** Let us assume that $\omega = 0$ and $\sigma, b, q, \psi \in \mathcal{S}(\mathbb{R}^d)$. Let us put $\Sigma = \Sigma_{\theta + \frac{\pi}{2}}$. Since $\sigma$ satisfies Hypothesis 1 and $\psi$ has generalized Blumenthal–Getoor index $s_1$, one sees by the product and chain rule, that $a(x, \xi) = \psi(\sigma^T(x)\xi) \in S_{1,0}^{s_1}(\mathbb{R}^d, \mathbb{R}^d)$. In particular,

(1) the mapping $(x, \xi) \mapsto a(x, \xi) := \psi(\sigma^T(x)\xi)$ belongs to $\mathcal{S}(\mathbb{R}^d \times \mathbb{R}^d)$;

---
[1]For the definition of $\operatorname{Hyp}_{\delta,\rho}^{m,m_0}(\mathbb{R}^d, \mathbb{R}^d)$, see Definition B.6, p. 20.

(2) for all multi-indices $\alpha$ and $\beta$ with $|\alpha|, |\beta| \leq k$ one has
$$\left|\partial_x^\beta \partial_\xi^\alpha \psi(\sigma^T(x)\xi)\right| \leq \langle|\xi|\rangle^{s_1-|\alpha|}, \quad |\xi|+|x| \geq 1;$$

(3) since $\sigma \asymp 1$, for all multi-indices $\alpha$ and $\beta$ with $|\alpha|, |\beta| \leq k$ one has
$$\left|\partial_x^\beta \partial_\xi^\alpha \psi(\sigma^T(x)\xi)\right| \leq |\xi|^{-|\alpha|}, \quad |\xi| \leq 1.$$

Furthermore, since $\psi$ has a generalized lower Blumenthal–Getoor index $s_1$, it follows that
$$\langle|\xi|\rangle^{s_1} \leq \left|\psi(\sigma^T(x)\xi)\right|, \quad \xi, x \in \mathbb{R}^d, \, |\xi| \geq R,$$
and, therefore, $\psi(\sigma^T(x)\xi) \in \text{Hyp}_{1,0}^{s,s_1}(\mathbb{R}^d, \mathbb{R}^d)$.

Let us denote by $a(x, D)$ the pseudo–differential operator induced by the symbol $a(x,\xi) = \psi(\sigma^T(x)\xi)$ and let us define the parameterized family of symbol $\{a(x,\xi,\lambda) := \psi(\sigma(x)^T\xi) + \lambda : \lambda \in \Sigma\}$. Observe, that $a(x, D, \lambda) := (\lambda + A)$ for all $\lambda \in \Sigma$. A short computation gives, that for all $\xi, x \in \mathbb{R}^d$ and $\lambda \in \Sigma$
$$\langle \lambda^{\frac{1}{s_1}} + |\xi|\rangle^{s_1} \lesssim |a(x,\xi,\lambda)|,$$
$$|a(x,\xi,\lambda)| \lesssim \langle \lambda^{\frac{1}{s_1}} + |\xi|\rangle^{s_1},$$
and
$$\left|\left(\partial_\xi^\alpha \partial_x^\beta a(x,\xi,\lambda)\right) 1/a(x,\xi,\lambda)\right| \lesssim \langle \lambda^{\frac{1}{s_1}} + |\xi|\rangle^{-|\alpha|},$$
for $\xi$ and $\lambda$ large enough.

It follows that $a(x,\xi,\lambda) \in \text{Hyp}L_{1,0,1/s_1}^{s_1,s_1}(\mathbb{R}^d, \mathbb{R}^d, \Sigma)$ (for the definition of $\text{Hyp}L$ see page 21) then by Theorem B.5 we know that the symbol $r(x,\xi,\lambda)$ of the resolvent $(A + \lambda I)^{-1} = R(\lambda, A)$ belongs to $\text{Hyp}L_{1,0,1/s_1}^{-s_1,-s_1}(\mathbb{R}^d, \mathbb{R}^d, \Sigma)$. Moreover, since the upper generalized Blumenthal–Getoor index of $q(\xi)$ is smaller or equal to $s_2$ and $|\partial_x^\alpha b(x)| \lesssim 1$ with $|\alpha| \leq k$, it follows that $\tilde{q}(x,\xi) := q(b^T(x)\xi) \in S_{1,0}^{s_2}(\mathbb{R}^d, \mathbb{R}^d)$. Therefore, by Theorem B.1, we know that the symbol $\tilde{\Theta}(x,\xi)$ of the composition of the two operators
$$\Theta(x, D) := B \, R(\lambda : A)$$
belongs to $S_{1,0}^{s_2-s_1}(\mathbb{R}^d, \mathbb{R}^d, \Sigma)$. Put $\Phi = \langle|\xi|+|x|\rangle$ and $\Psi(\xi, x) = 1$. Using the refined symbol class given in Definition B.9, one has $\tilde{q}(x,\xi) \in S(M_q, \phi, \psi)$ with $M_{\tilde{q}(x,\xi)} = \langle|\xi|+|x|\rangle^{s_2}$ and $r(x,\xi,\lambda) \in S(M_{r,\lambda}, \phi, \psi)$ with $M_{r,\lambda}(x,\xi) = \langle|\lambda|^{\frac{1}{s_1}} + |\xi|\rangle^{-s_1}$. In particular, by Theorem B.6 we get
$$\left|\tilde{\Theta}(x,\xi)\right| \lesssim \frac{\langle|\xi|+|x|\rangle^{s_2}}{\langle|\lambda|^{\frac{1}{s_1}} + |\xi|\rangle^{s_1}}, \quad x,\xi \in \mathbb{R}^d, \, \lambda \in \Sigma.$$

We first show that for all multi-indices $\alpha$ and $\beta$ with $|\alpha| \leq k$ and $|\beta| \leq k$ one has
$$\sup_{(x,\xi)\in\mathbb{R}^d\times\mathbb{R}^d\setminus\{0\}} \left|\partial_x^\beta \partial_\xi^\alpha \tilde{\Theta}(x,\xi)\right| |\xi|^{|\alpha|} < \infty.$$

First, to get an estimate independent of $\xi$ and $x$, we prove that there exists a constant $R > 0$ such that
$$\frac{|\xi|^{s_2}}{(1+|\lambda|^{\frac{1}{s_1}}+|\xi|)^{s_1}} \lesssim \lambda^{\frac{s_2}{s_1}-1}, \quad x, \xi \in \mathbb{R}^d, \lambda \in \Sigma_R := \Sigma \cap \{\lambda \in \mathbb{C} : |\lambda| \geq R\}.$$

In fact, by differentiating we infer that the function $f_\lambda$ defined by
$$f_\lambda(\xi) := \frac{\xi^{s_2}}{\lambda + \xi^{s_1}}, \quad \xi \in \mathbb{R},$$
attains it maximum at
$$|\xi_\lambda| = \left(\frac{s_2 \lambda}{s_1 - s_2}\right)^{\frac{1}{s_1}}.$$

Hence, there exists a constant $C = C(s_1, s_2) > 0$ such that
$$\left|\tilde{\Theta}(x,\xi)\right| \lesssim \lambda^{\frac{s_2}{s_1}-1}, \quad x, \xi \in \mathbb{R}^d, \lambda \in \Sigma.$$

This implies that $\tilde{\Theta}(x,\xi) \in S^0_{0,1}(\mathbb{R}^d, \mathbb{R}^d \setminus \{0\})$. Performing short and elementary calculations we can show that for all $\lambda \in \Sigma$
$$\sup_{(x,\xi) \in \mathbb{R}^d \times \mathbb{R}^d \setminus \{0\}} \left|\partial_x^\beta \partial_\xi^\alpha \tilde{\Theta}(x,\xi)\right| |\xi|^{|\alpha|} \lesssim \lambda^{\frac{s_2}{s_1}-1},$$
which altogether with Theorem B.4 implies (2.7).

In order to generalize the inequality to the Bessel Potential spaces $H_2^\rho(\mathbb{R}^d)$, where $\rho \in \mathbb{R}$, we notice first that for $\lambda \in \Sigma$ the symbol
$$\langle|\xi|\rangle^\rho \tilde{\Theta}(x,\xi) \langle|\xi|\rangle^{-\rho}$$
belongs also to $S^{s_2-s_1}_{0,0,1/s_1}(\mathbb{R}^d, \mathbb{R}^d)$ with norm
$$\|\langle|\xi|\rangle^\rho \tilde{\Theta}(x,\xi) \langle|\xi|\rangle^{-\rho}\|_{S^0_{k;0,0}} \leq \lambda^{\frac{s_2}{s_1}-1}.$$

Next, note that $(I+\Delta)^{\frac{\rho}{2}}$ is an isomorphism from $H_2^\tau(\mathbb{R}^d)$ to $H_2^{\tau+\rho}(\mathbb{R}^d)$. Thus for any $v \in H_2^\rho(\mathbb{R}^d)$ there exists a unique $u \in L^2(\mathbb{R}^d)$ such that $(I+\Delta)^{-\frac{\rho}{2}}u = v$. Now,
$$\left|(I+\Delta)^{\frac{\rho}{2}} B R(\lambda : A)(I+\Delta)^{-\frac{\rho}{2}}u\right|_{L^2} \leq C\, |\lambda|^{\varepsilon-1}\, |u|_{L^2}, \quad \forall u \in L^2(\mathbb{R}^d).$$

In addition, one has
$$\left|B R(\lambda : A)(I+\Delta)^{-\frac{\rho}{2}}u\right|_{H_2^\rho} = |B R(\lambda : A)v|_{H_2^\rho},$$
and
$$\left|(I+\Delta)^{\frac{\rho}{2}} B R(\lambda : A)(I+\Delta)^{-\frac{\rho}{2}}u\right|_{L^2} = \left|B R(\lambda : A)(1+\Delta)^{-\frac{\rho}{2}}u\right|_{H_2^\rho}.$$

By $|u|_{L^2} = |v|_{H_2^\rho}$ and Theorem B.4 we get
$$|B R(\lambda : A)u|_{H_2^\rho} \leq C\, |\lambda|^{\varepsilon-1}\, |u|_{H_2^\rho}, \quad \forall u \in H_2^\rho(\mathbb{R}^d).$$

Finally, Theorem A.1 gives the assertion for general $\rho \in \mathbb{R}$.

**Step 2:** If $\psi$ and $q$ belong to $C^k(\mathbb{R}^d)$, since $\mathcal{S}(\mathbb{R}^d)$ is dense in $C^k(\mathbb{R}^d)$, there exist sequences $\{\psi_n : n \in \mathbb{N}\}$, $\{\sigma_n : n \in \mathbb{N}\}$, $\{q_n : n \in \mathbb{N}\}$ and $\{b_n : n \in \mathbb{N}\} \subset \mathcal{S}(\mathbb{R}^d)$ converging to $\psi$, $\sigma$, $q$, and $b$, respectively. If $\psi, q \in C^k(\mathbb{R}^d \setminus \{0\})$, then we show that there exist $\{\psi_n : n \in \mathbb{N}\}$ and $\{q_n : n \in \mathbb{N}\} \subset \mathcal{S}(\mathbb{R}^d)$ such that for all multi-indices $\alpha$ and $\beta$ with $|\alpha| \leq k$ and $|\beta| \leq k$, the following limit holds

$$\sup_{(x,\xi) \in \mathbb{R}^d \times \mathbb{R}^d \setminus \{0\}} \lambda^{1-\varepsilon} |\xi|^\alpha \partial_\xi^\alpha \partial_x^\beta \left[ \frac{q_n(b_n^T(x)\xi)}{\lambda + \psi_n(\sigma_n(x)^T \xi)} - \frac{q(b^T(x)\xi)}{\lambda + \psi(\sigma(x)^T \xi)} \right] \to 0,$$

as $n \to \infty$. Let $h \in C^k(\mathbb{R}^d)$ be a function such that

- $h(x) = 0$, $x = 0$;
- $h(x) = 1$, $|x| \geq R$ for some $R > 0$;
- $|h|_{C_b^k(\mathbb{R}^d)} \leq 1$;
- $\lim_{|x| \to 0} |\partial_x^\alpha h(x)| = 0$ for all multi-index $\alpha$ with $|\alpha| \leq k$;
- $|\partial_x^\alpha h(x)| \leq |x|^{k-|\alpha|+1}$ for all multi-index $\alpha$ with $|\alpha| \leq k$.

Put

$$\tilde{\psi}_n(x) := h(xn)\psi(x), \quad x \in \mathbb{R}^d,$$

and

$$\tilde{q}_n(x) := h(xn)q(x), \quad x \in \mathbb{R}^d.$$

For each $n$ there exists sequences $\{\hat{\psi}_k^n : k \in \mathbb{N}\}, \{\hat{q}_k^n : k \in \mathbb{N}\} \subset \mathcal{S}(\mathbb{R}^d)$, such that $\hat{\psi}_k^n \to \tilde{\psi}_n$ in $C^k(\mathbb{R}^d)$ and $\hat{q}_k^n \to \tilde{q}_n$ in $C^k(\mathbb{R}^d)$. Put $\psi_n = \hat{\psi}_n^n$ and $q_n = \hat{q}_n^n$. Now, a straightforward calculations shows that for any multi-index $\alpha$ with $|\alpha| \leq k$ one has

$$\sup_{\xi \in \mathbb{R}^d \setminus \{0\}} |\xi|^{|\alpha|} \left| \partial_\xi^\alpha [\psi_n(\xi) - \psi(\xi)] \right| \to 0,$$

and

$$\sup_{\xi \in \mathbb{R}^d \setminus \{0\}} |\xi|^{|\alpha|} \left| \partial_\xi^\alpha [q_n(\xi) - q(\xi)] \right| \to 0,$$

for $n \to \infty$. Furthermore, one has

$$\sup_{\xi \in \mathbb{R}^d \setminus \{0\}} |\xi|^{|\alpha|} \left| \partial_\xi^\alpha \left[ \psi_n(\sigma_n^T(x)\xi) - \psi(\sigma^T(x)\xi) \right] \right| \to 0,$$

and

$$\sup_{\xi \in \mathbb{R}^d \setminus \{0\}} |\xi|^{|\alpha|} \left| \partial_\xi^\alpha \left[ q_n(b_n^T(x)\xi) - q(b^T(x)\xi) \right] \right| \to 0,$$

for $n \to \infty$. Let $\pi_n$ be the symbol of the composition $q_n \circ r_n$ where $q_n$ and $r_n$ are as above. Put

$$\pi_n(x, \xi, \lambda) = \frac{q_n(b_n^T(x)\xi)}{\lambda + \psi_n(\sigma_n^T(x)\xi)}.$$

Because of (B.1) and (B.2) one has

$$\sup_{\xi \in \mathbb{R}^d \setminus \{0\}} |\xi|^\alpha \left| \partial_\xi^\alpha [\pi_n(x,\xi) - \pi(x,\xi)] \right| \to 0$$

as $n \to \infty$, with
$$\pi(x,\xi,\lambda) = \frac{q(b^T(x)\xi)}{\lambda + \psi(\sigma^T(x)\xi)}.$$
Let us now replace in Step I $\psi$ and $q$ by the sequences $\{\psi_n : n \in \mathbb{N}\}$ and $\{q_n : n \in \mathbb{N}\}$ constructed above, and $\sigma$ and $b$ by sequences $\{\sigma_n : n \in \mathbb{N}\}$ and $\{b_n : n \in \mathbb{N}\}$ converging to $\sigma$ and $b$ in $C^k(\mathbb{R}^d)$. The Lebesgue Dominated Convergence Theorem gives the assertion.

To tackle the case where $\omega \neq 0$ it is sufficient to shift $\lambda$ and use a similar argument to the above. $\square$

## Appendix A. $C_0$–semigroups

For $\delta \in [0,\pi]$ we define $\Sigma_\delta := \{z \in \mathbb{C} \setminus \{0\} : |\arg(z)| < \delta\}$. An analytic $C_0$-semigroup is defined as follows (see also [22, Chapter 2.5]):

**Definition A.1.** *A degenerate $C_0$-semigroup $(T(t))_{t \geq 0}$ on a Banach space $U$ is called* analytic *in $\Sigma_\theta$ for some $\theta \in (0, \frac{\pi}{2})$ if*
  (1) *$T$ extends to an analytic function $T : \Sigma_{\theta+\frac{\pi}{2}} \to \mathcal{L}(U)$;*
  (2) *$\lim_{z \to 0; z \in \Sigma_\theta} T(z)x = T(0)x$ for all $x \in U$.*

Let
$$D(A) = \left\{ x \in U : \lim_{h \to 0} \frac{T(h)x - x}{h} \text{ exists} \right\}$$
and $A$ be an infinitesimal generator of a semigroup $T$, i.e.
$$Ax := \lim_{h \to 0} \frac{T(h)x - x}{h}, \quad x \in D(A).$$
Let $\rho(A)$ be the set of all complex numbers $\lambda$ for which $\lambda - A$ is invertible, i.e.
$$\rho(A) := \left\{ \lambda \in \mathbb{C} : (\lambda - A)^{-1} \text{ is a bounded operator} \right\}.$$
In order to characterize a $C_0$–semigroup, we introduce the following definition.

**Definition A.2.** *Let $U$ be a Banach space and let $A$ be the generator of a degenerate analytic $C_0$-semigroup on $U$. We say that $A$ is of type $(\omega, \theta, K)$, where $\omega \in \mathbb{R}$, $\theta \in (0, \frac{\pi}{2})$ and $K > 0$, if $\omega + \Sigma_{\frac{\pi}{2}+\theta} \subseteq \rho(A)$ and*
$$|\lambda - \omega| \|R(\lambda : A)\|_{L(U,U)} \leq K \quad \text{for all } \lambda \in \omega + \Sigma_{\frac{\pi}{2}+\theta}.$$

The theorem below gives some characterizations of analytic $C_0$-semigroups that we will use later on.

**Theorem A.1.** *Let $(T(t))_{t \geq 0}$ be a degenerate $C_0$-semigroup of type $(M, \omega)$ for some $M > 0$ and $\omega \in \mathbb{R}$. Let $A$ be the generator of $T$. Let $\omega' > \omega$. The following statements are equivalent:*
  (1) *$T$ is an analytic $C_0$-semigroup on $\Sigma_\theta$ for some $\theta \in (0, \frac{\pi}{2})$ and for every $\theta' < \theta$ there exists a constant $C_{1,\theta'}$ such that $\|e^{-\omega' z} T(z)\| \leq C_{1,\theta'}$ for all $z \in \Sigma_{\theta'}$.*

(2) *There exists a $\delta \in (0, \frac{\pi}{2})$ such that:*

$$\omega' + \Sigma_{\frac{\pi}{2}+\delta} \subset \rho(A),$$

*and for every $\delta' \in (0, \delta)$ there exists a constant $C_{2,\delta'} > 0$ such that:*

$$|\lambda - \omega'|\|R(\lambda : A)\| \leq C_{2,\delta'}, \quad \text{for all } \lambda \in \omega' + \Sigma_{\frac{\pi}{2}+\delta'}.$$

(3) *$T$ is differentiable for $t > 0$ and there exists a constant $C_3$ such that:*

$$t\|AT(t)\| \leq C_3 e^{\omega' t}, \quad \text{for all } t > 0.$$

The proof can be found in [22, Theorem 2.5.2] for exponentially stable analytic $C_0$-semigroups with bounded and invertible generator, and can be transferred to arbitrary analytic $C_0$– semigroup by the following observation. If $T$ is a $C_0$-semigroup of type $(M, \omega, K)$ and $A$ is the generator of $T$, then for any $\omega' > \omega$ the $C_0$-semigroup $(e^{-\omega' t}T(t))_{t\geq 0}$ is exponentially stable and the generator of this semigroup, $A - \omega' I_U$, is invertible.

For our purpose we need an estimate which is very similar to estimate (3) of Theorem A.1.

**Proposition A.1.** *Let $U$ be a Banach space. Let $A_0$ be the generator of a degenerate analytic $C_0$-semigroup $T$ on $U$ and let $B$ be a possible unbounded operator acting on $U$. Suppose $A_0$ is of type $(\omega, \theta, K)$ for some $\omega \in \mathbb{R}, \theta \in (0, \frac{\pi}{2})$ and $K > 0$. Suppose there exist an $\varepsilon \in [0, 1)$ and a constant $C(A_0, B)$ such that for all $\lambda \in \omega + \Sigma_{\frac{\pi}{2}+\theta}$ one has:*

$$\|R(\lambda : A_0)B\|_{\mathsf{L}(U,U)} \leq C(A_0, B)|\lambda - \omega|^{\varepsilon-1}.$$

*Then for all $t > 0$ one has*

$$\|T(s)B\|_{\mathsf{L}(U,U)} \leq \begin{cases} \frac{\Gamma(\varepsilon)}{\pi}[\sin\theta]^{-\varepsilon}C(A_0, B)\, e^{\omega s}s^{-\varepsilon} & \text{if } \varepsilon \neq 0, \\ \frac{1}{\pi}C(A_0, B)e^{\omega s}\int_0^\infty e^{-r\sin\theta}r^{-1}dr & \text{if } \varepsilon = 0. \end{cases}$$

*Proof of Proposition A.1.* First assume that $\omega = 0$. Let $\theta' \in (0, \theta)$, $\rho \in (0, \infty)$, and $\Gamma_{\theta',\rho} = \Gamma^{(1)}_{\theta',\rho} + \Gamma^{(2)}_{\theta',\rho} + \Gamma^{(3)}_{\theta',\rho}$, where $\Gamma^{(1)}_{\theta',\rho}$ and $\Gamma^{(2)}_{\theta',\rho}$ are the rays $re^{i(\frac{\pi}{2}+\theta')}$ and $re^{-i(\frac{\pi}{2}+\theta')}$, $\rho \leq r < \infty$, and $\Gamma^{(3)}_{\theta',\rho} = \rho^{-1}e^{i\phi}$, $\phi \in [-\frac{\pi}{2}-\theta', \frac{\pi}{2}+\theta']$. It follows from [22, Theorem 1.7.7] that for $s > 0$

$$T(s) = \tfrac{1}{2\pi i}\int_{\Gamma_{\theta'}(\rho)} e^{\lambda s}R(\lambda : A_0)d\lambda.$$

For any $x \in U$ one has

$$T(s)Bx = \tfrac{1}{2\pi i}\int_{\Gamma_{\theta'}} e^{\lambda s}R(\lambda : A_0)\,Bx\,d\lambda.$$

Therefore, for any $\varepsilon \in [0,1)$ and $s > 0$ one has the following chain of equalities/inequalities

$$\|T(s)\,B\|_{\mathsf{L}(U,U)} = \left\|\tfrac{1}{2\pi i}\int_{\Gamma_{\theta',\rho}} e^{\lambda s} R(\lambda:A_0)B\,d\lambda\right\|_{\mathsf{L}(U,U)}$$

$$\leq \frac{1}{2s\pi}\left\|\int_\rho^\infty e^{re^{-i(\frac{\pi}{2}+\theta')}} R(\tfrac{r}{s}e^{-i(\frac{\pi}{2}+\theta')}:A_0)\,B\,e^{i(\frac{\pi}{2}+\theta')}dr\right\|_{\mathsf{L}(U,U)}$$

$$+\frac{1}{2s\pi}\left\|\int_\rho^\infty e^{re^{i(\frac{\pi}{2}+\theta')}} R(\tfrac{r}{s}e^{i(\frac{\pi}{2}+\theta')}:A_0)\,B\,e^{-i(\frac{\pi}{2}+\theta')}dr\right\|_{\mathsf{L}(U,U)}$$

$$\frac{1}{2s\pi}\left\|\int_{-\frac{\pi}{2}-\theta'}^{\frac{\pi}{2}+\theta'} e^{\rho e^{i\alpha}} R(\tfrac{\rho}{s}e^{i\alpha}:A_0)\,B\,\rho e^{i\alpha}d\alpha\right\|_{\mathsf{L}(U,U)}$$

$$\leq \frac{1}{2s\pi}\int_\rho^\infty e^{-r\sin\theta'}\|R(\tfrac{r}{s}e^{-i(\frac{\pi}{2}+\theta')}:A_0)\,B\|_{\mathsf{L}(U,U)}dr$$

$$+\frac{1}{2s\pi}\int_\rho^\infty e^{-r\sin\theta'}\|R(\tfrac{r}{s}e^{i(\frac{\pi}{2}+\theta')}:A_0)\,B\|_{\mathsf{L}(U,U)}dr$$

$$+\frac{1}{2\pi s}\int_{-\frac{\pi}{2}-\theta'}^{\frac{\pi}{2}+\theta'} e^{\rho\cos\alpha}\|R(\tfrac{\rho}{s}e^{i\alpha}:A_0)\,B\|_{\mathsf{L}(U,U)}\rho\,d\alpha$$

$$\leq \frac{1}{s^\varepsilon \pi}C(A_0,B)\int_\rho^\infty e^{-r\sin\theta'} r^{\varepsilon-1}dr$$

$$+\rho^\varepsilon \frac{1}{2s^\varepsilon \pi}\int_{-\frac{\pi}{2}-\theta'}^{\frac{\pi}{2}+\theta'} e^{\rho\cos\alpha}d\alpha.$$

By letting $\rho \to 0$

$$\|T(s)\,B\|_{\mathsf{L}(U,U)} \leq \frac{1}{s^\varepsilon \pi}C(A_0,B)[\sin\theta']^{-\varepsilon}\int_0^\infty [r\sin\theta']^{\varepsilon-1}e^{-r\sin\theta'}\sin\theta'\,dr.$$

Therefore,

$$\|T(s)\,B\|_{\mathsf{L}(U,U)} \leq \begin{cases} \frac{\Gamma(\varepsilon)}{\pi}[\sin\theta']^{-\varepsilon}C(A_0,B)\,s^{-\varepsilon} & \text{if } \varepsilon \neq 0, \\ \frac{1}{\pi}C(A_0,B)\int_0^\infty e^{-r\sin\theta'}r^{-1}dr & \text{if } \varepsilon = 0. \end{cases}$$

□

## Appendix B. Symbol Classes and pseudo–differential operators

In this section we shortly introduce pseudo–differential operators and their symbols. In addition we introduce the definitions and theorems which are necessary to our purpose. For a detailed introduction on pseudo–differential operators and their symbols in the context of partial differential equations we recommend the books [21, 27, 30], in the context of Markov processes we recommend the books [17, 18, 19] or the survey [4].

In order to treat pseudo–differential operators different classes of symbols have been introduced. Here, we closely follow the definition of [27].

**Definition B.1.** *Let $U \subset \mathbb{R}^d$, $m \in \mathbb{R}$, and $\rho, \delta$ two real numbers such that $0 \leq \rho \leq 1$ and $0 \leq \delta \leq 1$. Let $S^m_{\rho,\delta}(U, \mathbb{R}^d)$ be the set of all functions $a : U \times \mathbb{R}^d \to \mathbb{C}$, where*

- *$a(x, \xi)$ is infinitely often differentiable, i.e. $a \in C^\infty(U \times \mathbb{R}^d)$;*
- *for any two multi-indices $\alpha$ and $\beta$ and any compact set $K \subset U$ there exists $C_{K,\alpha,\beta}$ such that*

$$\left|\partial_x^\alpha \partial_\xi^\beta a(x,\xi)\right| \leq C_{K,\alpha,\beta} \langle |\xi| + |x| \rangle^{m-\rho|\beta|+\delta|\alpha|}, \quad x \in U, \xi \in \mathbb{R}^d.$$

We call any function $a(x,\xi)$ belonging to $\cup_{m \in \mathbb{R}} S^m_{0,0}(\mathbb{R}^d, \mathbb{R}^d)$ a *symbol*. For many estimates, one does not need that the function is infinitely often differentiable. For this reason, one introduces also the following classes.

**Definition B.2.** *(compare [30, p. 28]) Let $U \subset \mathbb{R}^d$. Let $m \in \mathbb{R}$. Let $S^m_{k;\rho,\delta}(U, \mathbb{R}^d)$ be the set of all functions $a : U \times \mathbb{R}^d \to \mathbb{C}$, where*

- *$a(x, \xi)$ is $k$-times differentiable;*
- *and for any two multi-indices $\alpha$ and $\beta$ with $|\alpha| + |\beta| \leq k$, there exists a positive constant $C_{\alpha,\beta} > 0$ depending only on $\alpha$ and $\beta$ such that*

$$\left|\partial_x^\alpha \partial_\xi^\beta a(x,\xi)\right| \leq C_{\alpha,\beta} \langle |\xi| + |x| \rangle^{m-\rho|\beta|+\delta|\alpha|}, \quad x \in U, \xi \in \mathbb{R}^d.$$

Moreover, one can introduce a semi–norm in $S^m_{k;\rho,\delta}(U, \mathbb{R}^d)$ by

$$\|a\|_{k, S^m_{\delta,\rho}}$$
$$= \sup_{|\alpha|,|\beta| \leq k} \sup_{(x,\xi) \in U \times \mathbb{R}^d} \left|\partial_x^\alpha \partial_\xi^\beta a(x,\xi)\right| \langle |\xi|+|x|\rangle^{\rho|\beta|+\delta|\alpha|-m}, \quad a \in S^m_{k;\rho,\delta}(U, \mathbb{R}^d).$$

**Remark B.1.** *For $m_1 \geq m_2$ it follows that $S^{m_1}_{\rho,\delta}(\mathbb{R}^d, \mathbb{R}^d) \supseteq S^{m_2}_{\rho,\delta}(\mathbb{R}^d, \mathbb{R}^d)$ and $S^{m_1}_{k;\rho,\delta}(\mathbb{R}^d, \mathbb{R}^d) \supseteq S^{m_2}_{k;\rho,\delta}(\mathbb{R}^d, \mathbb{R}^d)$, $k \in \mathbb{N}$.*

**Definition B.3.** *(compare [30, p.28, Def. 4.2]) Let $a(x, \xi)$ be a symbol. The pseudo–differential operator $a(x, D)$ associated to $a(x, \xi)$ is defined by*

$$a(x, D)\, u(x) := \int_{\mathbb{R}^d} e^{ix^T \xi} a(x,\xi)\, \hat{u}(\xi)\, d\xi \quad u \in \mathcal{S}(\mathbb{R}^d).$$

The product of two pseudo–differential operators is again a pseudo–differential operator and can be characterized as follows.

**Theorem B.1.** *(compare [30, Theorem 6.1, p. 54], [21, Theorem 1.2.16, p. 31]) Let $a_1(x,\xi) \in S^{m_1}_{1,0}(\mathbb{R}^d, \mathbb{R}^d)$ and $a_2(x,\xi) \in S^{m_2}_{1,0}(\mathbb{R}^d, \mathbb{R}^d)$. Then the product $b(x, D) := a_1(x, D)\, a_2(x, D)$ is again a pseudo–differential operator such that*

$$b(x, \xi) \in S^{m_1+m_2}_{1,0}(\mathbb{R}^d, \mathbb{R}^d).$$

*Moreover,*

$$(\text{B.1}) \qquad b(x,\xi) \sim \sum_{|\alpha| \leq m_1 + m_2} \frac{(-i)^{|\alpha|}}{\alpha!} \left(\partial_\xi^\alpha a_1(x,\xi)\right) \left(\partial_x^\alpha a_2(x,\xi)\right).$$

*The equation* (B.1) *means that*

$$\text{(B.2)} \qquad b(x,\xi) - \sum_{|\alpha|\leq N} \frac{(-i)^{|\alpha|}}{\alpha!} \left(\partial_\xi^\alpha a_1(x,\xi)\right) \left(\partial_x^\alpha a_2(x,\xi)\right)$$

*belongs to* $S_{1,0}^{m_1+m_2-N}(\mathbb{R}^d \times \mathbb{R}^d)$ *for every positive integer* $N$.

**Definition B.4.** *Let $a(x,\xi)$ be a symbol in $S_{0,0}^m(\mathbb{R}^d, \mathbb{R}^d)$ and $a(x,D)$ the associated pseudo–differential operator. Suppose there exists a linear operator $a^*(x,D) : \mathcal{S}(\mathbb{R}^d) \to \mathcal{S}(\mathbb{R}^d)$ such that*

$$(a(x,D)f, g) = (f, a^*(x,D)g), \quad f,g \in \mathcal{S}(\mathbb{R}^d).$$

*Then we call $a^*(x,D)$ the formal adjoint operator of the operator $a(x,D)$.*

The existence of the formal adjoint is given by the following Theorem.

**Theorem B.2.** *(compare [15, Corollary 3.6, p. 803], [30, p. 62, Theorem 7.1]) For a symbol $a(x,\xi) \in S_{0,0}^m(\mathbb{R}^d, \mathbb{R}^d)$ there exists a symbol $a^*(x,\xi) \in S_{0,0}^m(\mathbb{R}^d, \mathbb{R}^d)$ such that the operator defined by*

$$a^*(x,D)u(x) := \int_{\mathbb{R}^d} e^{ix^T\xi} a^*(x,\xi)\, \hat{u}(\xi)\, d\xi, \quad u \in \mathcal{S}(\mathbb{R}^d),$$

*is the adjoint operator of $a(x,D)$. In addition, $a^*(x,\xi)$ has the following expansion*

$$\text{(B.3)} \qquad a^*(x,\xi) \sim \sum_\alpha \frac{(-i)^{|\alpha|}}{\alpha!} \left(\partial_x^\alpha \partial_\xi^\alpha \overline{a}\right)(x,\xi).$$

In the next theorem we give sufficient condition for a pseudo-differential operator $a(x,D)$ to be continuous.

**Theorem B.3.** *(compare [30, Theorem 9.7]) Let $a(x,\xi) \in S_{1,0}^0(\mathbb{R}^d, \mathbb{R}^d)$. Then, for any $1 < p < \infty$ the operator $a(x,D) : L^p(\mathbb{R}^d) \to L^p(\mathbb{R}^d)$ is a linear and bounded operator.*

In fact, analyzing the proof of Theorem 9.7 [30, p. 79] one can see that the condition of the differentiability at the origin can be relaxed. Here, it is important to mention that the proof relies on the Theorem 2.5 [14, p. 120] (see also Theorem 4.23 [1]), from which one can clearly see the extension of the Theorem 9.7 of [30] to symbols, whose derivatives have a singularity at $\{0\}$. Moreover, analyzing line by line of the proof of Theorem 9.7, one can give an estimate of the norm of the operator.

**Theorem B.4.** *Let $k > \frac{d}{2}$ and $a(x,\xi) \in C^k(\mathbb{R}^d \times \mathbb{R}^d \setminus \{0\})$. Moreover, let us assume that for any multi-indices $\alpha$ and $\beta$ with $|\alpha|, |\beta| \leq k$*

$$\sup_{|\alpha|,|\beta|\leq k} \sup_{\substack{(x,\xi)\in\mathbb{R}^d\times\mathbb{R}^d\setminus\{0\} \\ |x|+|\xi|\geq 1}} \left|\partial_\xi^\alpha \partial_x^\beta a(x,\xi)\right| < \infty,$$

*and*

$$|\partial_\xi^\alpha \partial_x^\beta a(x,\xi)| \lesssim |\xi|^{-|\alpha|}, \quad \xi \neq 0.$$

*Then the corresponding operator $a(x,D)$ is bounded on $L^2(\mathbb{R}^d)$ with the uniform estimate*

$$\|a(x,D)\|_{L(L^2(\mathbb{R}^d))} \lesssim \sup_{|\alpha|,|\beta| \leq k} \sup_{(x,\xi) \in \mathbb{R}^d \times \mathbb{R}^d \setminus \{0\}} \left| |\xi|^{|\alpha|} \partial_x^\beta \partial_\xi^\alpha a(x,\xi) \right|.$$

To investigate the inverse of a pseudo–differential operator one should introduce the set of elliptic and hypoelliptic symbols.

**Definition B.5.** *(compare [21, p. 35]) A symbol $a \in S_{\rho,\delta}^m$ is called globally elliptic in the class $S_{\rho,\delta}^m(\mathbb{R}^d, \mathbb{R}^d)$, if for some $R > 0$,*

$$\langle |\xi| \rangle^m \lesssim |a(x,\xi)|, \quad |x| + |\xi| \geq R.$$

**Definition B.6.** *(compare [21, p. 35]) Let $m, m_0, \rho, \delta$ be real numbers with $0 \leq \delta < \rho \leq 1$. The class $\mathrm{Hyp}_{\rho,\delta}^{m,m_0}(U \times \mathbb{R}^d)$ consists of all functions $a(x,\xi)$ such that*
- $a(x,\xi) \in C^\infty(U \times \mathbb{R}^d)$;
- *there exists some $R > 0$ such that*

$$\langle |\xi| \rangle^{m_0} \lesssim |a(x,\xi)|, \quad |x| + |\xi| \geq R.$$

*and for an arbitrary multi-indices $\alpha$ and $\beta$ and for any compact set $K \subset U$ there exists a constants $C_{\alpha,\beta,K}$ with*

$$\left| \partial_\xi^\alpha \partial_x^\beta a(x,\xi) \right| \leq C_{K,\alpha,\beta} \langle |\xi| + |x| \rangle^{m-\rho|\alpha|+\delta|\beta|}.$$

*for $x \in K$, $\xi \in \mathbb{R}^d$.*

**Remark B.2.** *Let $a(x,\xi) \in \mathrm{Hyp}_{\rho,\delta}^{m,m_0}(U \times \mathbb{R}^d)$ be a symbol and $a^*(x,\xi)$ the symbol of the formal adjoint operator $a^*(x,D)$. Then, one can see from the expansion in (B.3), that if $a(x,\xi) \in \mathrm{Hyp}_{\rho,\delta}^{m,m_0}(U \times \mathbb{R}^d)$ then $a^*(x,\xi) \in \mathrm{Hyp}_{\rho,\delta}^{m,m_0}(U \times \mathbb{R}^d)$.*

Lemma 1.3.5 [21, p. 36] gives under which conditions a symbol belonging to $\mathrm{Hyp}_{\rho,\delta}^{m,m_0}(U \times \mathbb{R}^d)$ is invertible. However, our aim is to characterize the symbol of a resolvent of an operator $a(x,D)$. Then we define a subclass of hypoelliptic symbols and state a theorem giving sufficient condition for the existence of the symbol of the resolvent.

**Definition B.7.** *(compare [27, Definition 9.1, p. 77]) Let $m, \rho, \delta, \gamma$ be real numbers with $0 \leq \delta < \rho \leq 1$, $0 < \gamma < \infty$. The class $S_{\rho,\delta;d}^m(U \times \mathbb{R}^d, \Lambda)$ consists of all functions $a(x,\xi,\lambda) : U \times \mathbb{R}^d \times \Lambda \to \mathbb{C}$ such that*
- $a(x,\xi,\lambda_0) \in C^\infty(U \times \mathbb{R}^d)$ *for every fixed $\lambda_0 \in \Lambda$;*
- *For arbitrary multi-indices $\alpha$ and $\beta$ and for any compact set $K \subset U$ there exists a constants $C_{K,\alpha,\beta}$ such that*

$$\left| \partial_\xi^\alpha \partial_x^\beta a(x,\xi,\lambda) \right| \leq C_{K,\alpha,\beta} \langle |\xi| + |\lambda|^{\frac{1}{\gamma}} \rangle^{m-\rho|\alpha|+\delta|\beta|}.$$

for $x \in K$, $\xi \in \mathbb{R}^d$, $\lambda \in \Lambda$.

Since the resolvent can be viewed as a parameterized family of symbols, we introduce the following definition.

**Definition B.8.** *(compare [27, p. 78]) Let $m, m_0, \rho, \delta, \gamma$ be real numbers with $0 \leq \delta < \rho \leq 1$, $0 < \gamma < \infty$. The class $Hyp_{\rho,\delta;\gamma}^{m,m_0}(U \times \mathbb{R}^d, \Lambda)$ consists of all functions $a(x, \xi, \lambda) : U \times \mathbb{R}^d \times \Lambda \to \mathbb{C}$ such that*

- *$a(x, \xi, \lambda_0) \in C^\infty(U \times \mathbb{R}^d)$ for every fixed $\lambda_0 \in \Lambda$;*
- *For any compact set $K \subset U$ there exists two constants $C_K$ and $\tilde{C}_K$ such that*

$$C_K \langle |\xi| + |\lambda|^{\frac{1}{\gamma}} \rangle^{m_0} \leq |a(x,\xi,\lambda)| \leq \tilde{C}_K \langle |\xi| + |\lambda|^{\frac{1}{\gamma}} \rangle^m.$$

*for $x \in K$, $\xi \in \mathbb{R}^d$, $\lambda \in \Lambda$, $|\xi| + |\lambda| \geq R$.*
- *For arbitrary multi-indices $\alpha$ and $\beta$ and for any compact set $K \subset U$ there exists a constants $C_{K,\alpha,\beta}$ such that*

$$\left| \left( \partial_\xi^\alpha \partial_x^\beta a(x,\xi,\lambda) \right) a^{-1}(x,\xi,\lambda) \right| \leq C_{K,\alpha,\beta} \langle |\xi| + |\lambda|^{\frac{1}{\gamma}} \rangle^{m-\rho|\alpha|+\delta|\beta|},$$

*for $x \in K$, $\xi \in \mathbb{R}^d$, $\lambda \in \Lambda$.*

**Remark B.3.** *Let $a(x, \xi) \in Hyp_{\rho,\delta}^{m,m_0}(\mathbb{R}^d \times \mathbb{R}^d)$ of type $(\theta, \omega)$. Then, it is easy to see that $a(x,\xi) + \lambda \in Hyp_{\rho,\delta;\frac{1}{m_1}}^{m,m_0}(\mathbb{R}^d \times \mathbb{R}^d \times \Lambda)$ with $m_1 = \max(m_0, m)$ and $\Lambda = \omega + \Sigma_{\frac{\theta}{2}}$.*

One can classify the inverse of the each member of $\{a(x, \xi, \lambda) : \lambda \in \Lambda\}$, but one has to introduce the set of properly supported operators. Let $a(x, \xi)$ be a symbol with kernel $K_a$ and let $\text{supp}(K_a)$ denote the support of $K_a$ (the smallest closed subset $Z \subset U \times U$ such that $K_a |_{(U \times U) \setminus Z} = 0$). Consider the canonical projections $\Pi_1, \Pi_2 : \text{supp}(K_a) \to U$, obtained by restricting the corresponding projections of the direct product $U \times U$. Recall that a continuous map $f : M \to N$ between topological spaces $M$ and $N$ is called proper if for any compact $K \subset N$ the inverse image $f^{-1}(K)$ is a compact in $M$. A symbol $a(x, \xi)$ is called properly supported if both projections $\Pi_1, \Pi_2 : \text{supp}(K_a) \to U$ are proper maps. For more details see [27, p. 18]. We will denote by $HypL_{\rho,\delta;\gamma}^{m,m_0}(U \times \mathbb{R}^d, \Lambda)$ the class of symbols $a \in Hyp_{\rho,\delta;\gamma}^{m,m_0}(U \times \mathbb{R}^d, \Lambda)$ being properly supported.

**Theorem B.5.** *(compare [27, Theorem 9.2, p. 85]) Let $a(x,\xi,\lambda) \in HypL_{\rho,\delta;d}^{m,m_0}(U \times \mathbb{R}^d, \Lambda)$. Then, there exists a $R > 0$ such that for all $|\lambda| \geq R$ the operator $a(x, D, \lambda)$ is invertible. In particular, there exists a symbol $a_{inv}(x, \xi, \lambda)$ such that*

$$a(x,D,\lambda)^{-1} u = \int_{\mathbb{R}^d} e^{ix^T \xi} a_{inv}(x,\xi,\lambda) \, \hat{u}(\xi) \, d\xi \quad u \in \mathcal{S}(\mathbb{R}^d),$$

*and $a_{inv}(x, D, \lambda)$ belongs to $HypL_{\rho,\delta;d}^{-m_0,-m}(U \times \mathbb{R}^d, \Lambda)$, where $\Lambda_R := \Lambda \cap \{\lambda \in \mathbb{C}, |\lambda| \geq R\}$.*

In order to deal with operators depending on a parameter, one can treat a family of symbols by refining the definition of symbol classes see for instance, [21, p. 19]. For this aim, we introduce some class of functions. A positive continuous function $\Phi : \mathbb{R}^d \times \mathbb{R}^d \to \mathbb{R}^d$ is called *sub–linear weight* function if,
$$1 \leq \Phi(x,\xi) \lesssim 1 + |x| + |\xi|, \quad \text{for } x, \xi \in \mathbb{R}^d.$$
It is called a *temperate weight*, if for some $s > 0$
$$\Phi(x+y, \xi+\eta) \lesssim \Phi(x,\xi) \left(1 + |y| + |\eta|\right)^s, \quad \text{for } x, y, \xi, \eta \in \mathbb{R}^d.$$

**Definition B.9.** *Let $\Phi, \Psi, M : \mathbb{R}^d \times \mathbb{R}^d \to \mathbb{R}^d$ be temperate weights such that $\Phi, \Psi$ are sub–linear. We denote by $S(M; \Phi, \Psi)$ the space of all smooth functions $a : \mathbb{R}^d \times \mathbb{R}^d \to \mathbb{R}^d$, such that for every $\alpha, \beta \in \mathbb{N}^d$ one has*
$$|\partial_\xi^\alpha \partial_x^\alpha a(x,\xi)| \lesssim M(x,\xi) \Psi(x,\xi)^{-|\alpha|} \Phi(x,\xi)^{-|\beta|}, \quad \text{for } x, y, \xi, \eta \in \mathbb{R}^d.$$

Now, one can choose the target weight $M$ in such a way that it depends on a parameter, say $\lambda$. Now the multiplication theorem for the composition operator can be restated as follows.

**Theorem B.6.** *([21, Theorem 1.2.16, p. 31]) Let $\Phi, \Psi, M : \mathbb{R}^d \times \mathbb{R}^d \to \mathbb{R}^d$ be temperate weights such that $\Phi, \Psi$ are sub–linear. Let $a_1(x, \xi) \in S(M_1; \Phi, \Psi)$ and $a_2(x, \xi) \in S(M_2; \Phi, \Psi)$. Then the composition $b(x, D) := a_1(x, D) a_2(x, D)$ is again a pseudo–differential operator such that*
$$b(x, \xi) \in S(M_1 M_2; \Phi, \Psi).$$

## Acknowledgements

This work was supported by the Austrian Science foundation (FWF), Project number P23591-N12. We would like to thank the anonymous referee for his insightful comments and remarks which improved the manuscript. We also thank him/her for providing us with Example 2.4.

ANALYTIC PROPERTIES OF MARKOV SEMIGROUP GENERATED BY SDES 23[8] S. Ethier and T. Kurtz. *Markov processes: Characterization and convergence.* Wiley Series in Probability and Mathematical Statistics: Probability and Mathematical Statistics. John Wiley & Sons, Inc., New York, (1986).

[9] G. Da Prato. *An introduction to infinite-dimensional analysis*, Springer-Verlag, Berlin, (2006).

[10] B. P. W. Fernando and E. Hausenblas. Nonlinear filtering with correlated Lévy noise characterized by copulas, submitted for publication.

[11] K. Glau. Sobolev index: a classification of Lévy processes via their symbols, arXiv:1203.0866, (2012).

[12] A. Gomilko and Y. Tomilov. On subordination of holomorphic functions, arXiv:1408.1417, (2014).

[13] M. Haase. *The functional calculus for sectorial operators*, Volume 169 of *Operator Theory: Advances and Applications*, 169. Birkhäuser Verlag, Basel, (2006).

[14] L. Hörmander. Estimates for translation invariant operators in $L^p$ spaces, *Acta Math.*, 104,93–140 (1960).

[15] W. Hoh. A symbolic calculus for pseudo-differential operators generating Feller semigroups, Osaka J. Math., 35, 798–820 (1998).

[16] N. Jacob and R. L. Schilling. Fractional derivatives, non-symmetric and time-dependent Dirichlet forms and the drift form, *Z. Anal. Anwendungen*,(19)3, 801830 (2000).

[17] N. Jacob. *Pseudo differential operators and Markov processes*, Vol. I, Fourier analysis and semigroups, Imperial College Press, London (2000).

[18] N. Jacob. *Pseudo differential operators and Markov processes*, Vol. II, Generators and their potential theory, Imperial College Press, London (2002).

[19] N. Jacob. *Pseudo differential operators and Markov processes*, Vol. III, Markov processes and applications, Imperial College Press, London (2005).

[20] U. Küchler and S. Tappe. *Tempered stable distributions and processes*, Stochastic Process. Appl., 123(12), 4256–4293 (2013).

[21] F. Nicola and L. Rodino. *Global pseudo-differential calculus on Euclidean spaces: Pseudo-Differential Operators,* Theory and Applications, 4. Birkhäuser Verlag, Basel, (2010).

[22] A. Pazy. *Semigroups of Linear Operators and Applications to Partial Differential Equations*, Applied Mathematical Sciences, 44, Springer-Verlag, New York, (1983).

[23] K. Sato. *Lévy processes and infinitely divisible distributions,* Cambridge Studies in Advanced Mathematics, 68. Cambridge University Press, Cambridge, (1999).

[24] R. Schilling. *Growth and Hölder conditions for the sample paths of Feller processes.* Probab. Theory Rel. Fields 112 , 565-611 (1998).

[25] R. Schilling, and A. Schnurr. *The symbol associated with the solution of a stochastic differential equation*, Electron. J. Probab., 15, 1369-1393 (2010).

[26] W. Schoutens. *Meixner Processes in Finance*, http://alexandria.tue.nl/repository/books/548458.pdf.

[27] M. A. Shubin. *Pseudodifferential operators and spectral theory*, Springer-Verlag, Berlin, (2001).

[28] E. Stein. *Harmonic analysis: real-variable methods, orthogonality, and oscillatory integrals*, With the assistance of Timothy S. Murphy. Princeton Mathematical Series 43, Monographs in Harmonic Analysis III, Princeton University Press, Princeton, NJ, (1993).

[29] D. Haroske and H. Triebel. *Distributions, Sobolev spaces, elliptic equations*, Zürich: European Mathematical Society, (2008).

[30] M. W. Wong. *An introduction to pseudo-differential operators,* Second ed., World Scientific Publishing Co. Inc., River Edge, NJ, (1999).

*E-mail address*, P. Fernando: `fernando.bandhisattambige@unileoben.ac.at`

*E-mail address*, E. Hausenblas: `erika.hausenblas@unileoben.ac.at`

*E-mail address*, P. Razafimandimby: `paul.razafimandimby@unileoben.ac.at`

Lehrstuhl für Angewandte Mathematik,, Montanuniversität Leoben, Franz Josef Strasse 18, 8700 Leoben, Austria